\documentclass[12pt,a4paper,reqno]{amsart}
\usepackage{geometry}
\usepackage[T1]{fontenc}

\usepackage{enumerate, mathabx, mathtools}
\usepackage{amsthm, amsmath, amsfonts, amssymb}

\usepackage{cite}
\usepackage{amssymb, amsmath, amsthm, enumerate, mathtools}
\usepackage{color}
\usepackage[dvipsnames]{xcolor}
\usepackage{graphicx}
\usepackage{tikz}

\usepackage{subcaption} 
\usepackage{caption}

\usepackage[thinlines]{easytable}

\usepackage{mathabx} 
\usepackage{esint} 
\usepackage{bm}
\usepackage{tikz}
\usepackage{tikz-3dplot}
\usepackage[normalem]{ulem}

\usetikzlibrary{calc,patterns,angles,quotes}

\usepackage{hyperref}
\usepackage{cite}

\usepackage{etoolbox}

\preto{\section}{}

\makeatletter
\@addtoreset{theorem}{section}
\makeatother

\usepackage{url}

\pdfoutput=1
\usepackage[latin1]{inputenc}
\usepackage{ae,aecompl,amsbsy,amssymb,amsmath,amsthm,
eurosym,amsfonts,epsfig,graphicx,graphics,verbatim,enumerate,color}
\usepackage{hyperref}
\hypersetup{colorlinks=true,linkcolor=blue}
\usepackage{tikz}
\usetikzlibrary{arrows, calc, positioning, decorations.markings}
\tikzset{
    between/.style args={#1 and #2}{
         at = ($(#1)!0.5!(#2)$)
    }
}

\let\oldaddcontentsline\addcontentsline

\newcommand{\starttocentries}{\let\addcontentsline\oldaddcontentsline}

\usepackage{enumitem}

\theoremstyle{plain}
\newtheorem{theorem}{Theorem}[section]

\theoremstyle{definition}

\theoremstyle{remark}
\newtheorem*{remark}{Remark}

\begin{document}
\date{}
\title[Littlewood--Paley theory and quadratic expressions]{Littlewood, Paley and Almost-Orthogonality: \\
a theory well ahead of its time} 
\author{Anthony Carbery}

\begin{abstract}
The classic paper J.\ E.\ Littlewood, R.\ E.\ A.\ C.\ Paley, {\em Theorems on Fourier Series and Power Series}. J. Lond. Math. Soc. (1), {\bf 6} (1931), 230--33, marked the birth of Littlewood--Paley theory.  We discuss this paper and its impact from a historical perspective. We include an outline of the results in the paper and their subsequent significance in relation to developments over the last century, and set them into the context of the current state of the art in harmonic analysis and beyond.
\end{abstract}

\address{Anthony Carbery, School of Mathematics and Maxwell Institute for Mathematical Sciences, University of Edinburgh, James Clerk Maxwell Building, Peter Guthrie Tait Rd, Kings Buildings, Edinburgh EH9 3FD, Scotland}
\email{A.Carbery@ed.ac.uk}


\maketitle
\setcounter{tocdepth}{2}
%
%
\section{Introduction}\label{sec:intro}

Euclidean geometry, and more generally the study of euclidean spaces, is rooted in the notion of orthogonality.
The same holds for the study of Hilbert spaces, the infinite-dimensional generalisations of euclidean spaces, of which $L^2$ is the prototypical example. The fundamental and beautiful spectral theorem for bounded linear self-adjoint operators on a Hilbert space is entirely founded upon the concept of orthogonality. However, many problems occurring in analysis, mathematical physics, geometry, number theory and PDE, {\em inter alia}, require consideration of $L^p$ spaces, where $p \in (1, \infty)$ is an index different from $2$. Is there a practical and useful substitute for orthogonality in such a setting, which might for example lead to results with flavours similar to those that the spectral theorem yields in the case $p=2$?

The foundational paper of Littlewood and Paley, \cite{LP1}, begins to address this question. An entire mathematical culture and body of thought, now going by the name Littlewood--Paley theory, can be traced back to to it. In the present article we review this paper, and sketch some of its ramifications over the last nearly 100 years. In a nutshell, the paper introduced the notion that ideas and techniques based on orthogonality can be developed and deployed in settings where one has no right to believe that they should even exist.

For the sake of concreteness, let us initially assume that we are working with $L^p$ spaces defined on a locally compact abelian group. For Littlewood and Paley, this group was the circle group $\mathbb{T}$, but we shall also consider euclidean space $\mathbb{R}^n$. The duals of these groups are $\mathbb{Z}$ and $\mathbb{R}^n$ respectively, and in both cases we have a Fourier transform $f \mapsto \widehat{f}$ which is an isometric isomorphism of $L^2$ of the group onto $L^2$ of  its dual.\footnote{Starting with $\mathbb{Z}$ and its dual group $\mathbb{T}$ leads to considerations of a different, fundamentally number-theoretic nature, culminating in the Hardy--Littlewood circle method, into which we do not enter here.}

As we have observed, orthogonality is baked into the Hilbert space $L^2$, and to operators on $L^2$ spaces. Two functions $f$ and $g$ in $L^2$ are orthogonal if and only if 
\[ \langle f, g \rangle = \int f \overline{g} =0
\]
and, in the setting of locally compact abelian groups, this is equivalent to
\[ \langle \widehat{f}, \widehat{g} \rangle = \int \widehat{f} \overline{\widehat{g}} =0.
\]
So we see that if $\{f_j\}$ are such that their Fourier transforms $\{\widehat{f}_j\}$ have pairwise disjoint supports, then, by Pythagoras' theorem, 
\[ 
\| \sum_j f_j \|_2^2 = \sum_j \|f_j\|_2^2.
\]
Put slightly differently, if $\{E_j\}$ forms a partition of Fourier space (the dual group), and if we define Fourier multiplier operators $Q_j$ by 
\[ \widehat{Q_j f}(\xi) := \chi_{E_j}(\xi) \widehat{f}(\xi),
\]
then we have the identity
\begin{equation}\label{eq:id} 
\|f\|_2^2 = \| \sum_j Q_j f\|_2^2 = \sum_j \|Q_jf\|_2^2 = \| \big(\sum_j |Q_j f|^2\big)^{1/2}\|_2^2.
\end{equation}

The philosophy which emerges from this perspective is that one may take a function $f$ or an operator $T$, (which we shall assume for the moment is a convolution operator),  
and decompose it into orthogonal pieces $Q_jf$, or
$Q_jT = TQ_j$, work with each piece separately, and then use Pythagoras' theorem to re-assemble the pieces.

When confronted instead with a problem whose correct expression is not in the language of $L^2$-spaces, but instead of that of $L^p$ spaces, where $1 < p < \infty$ and $p \neq 2$, we seem doomed to failure if we entertain the idea of trying to analyse it using techniques of orthogonality. The genius of Littlewood and Paley,
as manifested in this paper \cite{LP1}, is encapsulated in the notion that for $f \in L^p$, by considering the {\em quadratic expression} 
\[\mathcal{Q}f (x)= \big(\sum_j |Q_j f(x)|^2\big)^{1/2}
\]
which appears in \eqref{eq:id}, and its quantitative relationship to $f$, we may indeed bring ideas of orthogonality to bear on problems which are essentially $L^p$ in nature.

The aim of this article is to give a taste of the story of these and related quadratic expressions, or square functions, the $L^p$ estimates they satisfy, and a hint of their application and influence in the wider context of almost-orthogonality, over the last nearly one hundred years since Littlewood--Paley theory first emerged in \cite{LP1}.

{\em Disclaimer and acknowledgements.} What follows is a 
personal, subjective and in places rather impressionistic account, which reflects the tastes, interests and experience of its author. There is no claim (explicit or implicit) to comprehensiveness in either the text or the references. In particular there will be many authors whose work relies upon and is intimately bound up with Littlewood--Paley theory, and which indeed advances it in significant ways, but who are not mentioned here; I can only apologise in advance for such omissions. We make no attempt to address the endpoint cases $p=1$ and $p=\infty$, still less the cases $0 < p <1$. (These too have a very rich history, bound up with Hardy spaces $H^p$ in their classical and real-variable incarnations, and the John--Nirenberg space BMO of functions of bounded mean oscillation which is the dual space of the real Hardy space $H^1$. See \cite{SteinBig}.) For a complementary perspective focusing on Zygmund's work on square functions see \cite{SteinZygmund}; that article gives an authoritative and thorough account of matters up to the early 1980s. See also \cite{CarberyCZsurvey} for an account of related developments up to 1993. The author would like to thank Jonathan Hickman for comments on an earlier version of the paper, and Joseph Feneuil, Linhan Li and Michele Villa for guidance on Section~\ref{sec:GMT}. The comments of the referee have led to many presentational improvements. All errors, omissions and misrepresentations remaining are the sole responsibility of the author.

\section{What's in the Littlewood--Paley paper (and what's not)}\label{sec:secondsection}
In this paper \cite{LP1}, the programme of bringing orthogonality techniques to bear on $L^p$ for $p \neq 2$ is addressed in the specific setting of $L^p(\mathbb{T})$ for $1 < p < \infty$. Here, $\mathbb{T} = \mathbb{R}/\mathbb{Z}$ is the unit torus or circle group which we can identify with $(-1/2, 1/2]$, and the Fourier transform $f \mapsto \widehat{f}$ is the Fourier coefficient mapping from $L^2(\mathbb{T})$ to $\ell^2(\mathbb{Z})$ given by 
\[ \widehat{f}(k) = \int_\mathbb{T} f(x) e^{-2 \pi i k x} {\rm d} x.
\]
We shall specifically take the partitioning sets $E_j \subseteq \mathbb{Z}$ referred to above to comprise the {\em standard dyadic partition} of the integers; that is, $E_j := \mathbb{Z} \cap \big([2^{j-1}, 2^{j}) \cup (-2^{j}, -2^{j-1}]\big)$ for $j \geq 1$, together with $E_{0} = \{0\}$. The set $E_j$ is referred to as the $j$'th dyadic block of integers. Thus $Q_{0}f = \int_\mathbb{T} f$, while for $j \geq 1$
\[ Q_j f (x) = \sum_{2^{j-1} \leq |k| < 2^{j}} \widehat{f}(k) e^{2 \pi i kx }. 
\]
Note that the $Q_j$ are projection operators onto mutually orthogonal subspaces of $L^2(\mathbb{T})$.

\medskip
\noindent
The main result announced in the paper is the following:
\begin{theorem}\label{thm:LP}
For each $ 1 < p < \infty$, there are constants $A_p$ and $B_p$, such that for all $f \in L^p$,
\begin{equation}\label{eq:LP}
A_p  \|f\|_p \leq \| \big(\sum_{j \geq 0} |Q_j f|^2\big)^{1/2}\|_p \leq  B_p  \|f\|_p.
\end{equation}
\end{theorem}
Thus, at the loss of a constant depending only on $p$, one may replace any instance of 
$\|f\|_p$ by the $L^p$-norm of the quadratic expression $\mathcal{Q}f$. This represents the first and essential step in the programme outlined above. Because of its links with identity \eqref{eq:id}, Theorem~\ref{thm:LP} is often referred to as an ``almost-orthogonality" result. The second inequality in \eqref{eq:LP} is the key result; once one has this for a given $1 < p < \infty$ the first inequality follows for the dual index $p'$ (here, as always, $\frac{1}{p} + \frac{1}{p'} =1$) from \eqref{eq:id} by polarisation and duality.

\begin{remark}\label{rem:lac}
Note that this result immediately implies that for a {\em lacunary} trigonometric series $\sum_{j \geq 0}  a_j e^{2 \pi i 2^j x}$, its $L^p$ norms for $1 < p < \infty$ are all comparable to $(\sum_{j \geq 0} |a_j|^2)^{1/2}$, since each dyadic block $E_j$ then contains at most one non-zero Fourier frequency. The ``digitalised" version of this consequence -- in which the term $\sin(2^{j+1} \pi x)$ implicitly appearing in a lacunary trigonometric series is replaced by its signum -- goes back a little further to Khintchine \cite{MR1544623}. Results such as this capture the concept of ``square-root cancellation" which serves as an intuitional lynchpin across modern harmonic analysis and beyond.
\end{remark} 

As a pretty much immediate consequence of Theorem~\ref{thm:LP}, Littlewood and Paley derive:
\begin{theorem}\label{thm:Marc_primordial}
For each $ 1 < p < \infty$, there is a constant $C_p$ such that if the numerical sequence $\{r_j\}$ satisfies $|r_j| \leq 1$ for all $j$, then, for all $f \in L^p$,
\[ 
\| \sum_{j \geq 0}  r_j Q_j f \|_p \leq  C_p  \|f\|_p.
\]
\end{theorem}
Interestingly, in \cite{LP1} this result is stated under the stronger hypothesis that $|r_j| = 1$ for all $j$, and, in that case, the stronger conclusion that $ \|\sum_{j \geq 0}  r_j Q_j f \|_p  \approx \|f\|_p$ is not mentioned. Here as below, the symbol $\approx$ indicates that the ratio of the left and right hand sides is bounded above and below by constants depending only on $p$ (and the underlying dimension when that enters the story later).

The proof of Theorem~\ref{thm:Marc_primordial} relies upon applications of both inequalities in \eqref{eq:LP}, and can be summarised by
the chain of inequalities
\[
\| \sum_{j \geq 0}  r_j Q_j f \|_p 
\leq A_p^{-1} \| \big( \sum_{k \geq 0}  |Q_k (\sum_{j \geq 0} r_j Q_j f)|^2\big)^{1/2}\|_p \leq A_p^{-1} \| \big( \sum_{k \geq 0}  |Q_k f|^2\big)^{1/2}\|_p 
\leq A_p^{-1} B_p \|f\|_p,
\]
holding since the projection operators $Q_k$ satisfy $Q_k^2 = Q_k$ and $Q_kQ_j =0$ for $j \neq k$.

Theorem~\ref{thm:Marc_primordial} is a primordial version of subsequently established {\em Fourier multiplier theorems} going by the names of Marcinkiewicz\footnote{In this case, its ultimate formulation is still a matter of current study, see for example \cite{bakasetal}.} and H\"ormander--Mikhlin, extensive discussion of both of which (in the euclidean case) can be found in \cite{Stein}. 
In fact, the proofs of these more general multiplier theorems are not much more difficult than that of Theorem~\ref{thm:Marc_primordial}, as we now demonstrate. 

Indeed, let $T_j$ be a {\em Fourier multiplier operator} (defined on functions on $\mathbb{T}$) given by
\[ \widehat{T_j f}(k) = m_j(k) \widehat{f}(k),
\]
where $m_j$ is supported on $E_j$. The basic premise is that we shall suppose that the operators $T_j$ are uniformly (in $j$) bounded on $L^p$ for $1 < p < \infty$, and we wish to conclude that $\sum_{j\geq 0}T_j$ is also $L^p$ bounded in the same range of $p$. Such a conclusion would be reminiscent of what can be obtained from the spectral theorem when $p=2$.

In fact we shall impose the slightly stronger uniform weighted hypothesis
\begin{equation}\label{eq:unif}
\sup_j \int_{\mathbb{T}} |T_jf|^2 w \leq \int_\mathbb{T} |f|^2 \mathcal{M}w
\end{equation}
where, for each $1 < r < \infty$, $w \mapsto \mathcal{M}w$ is some (possibly non-linear) non-negative operator which is bounded on $L^r(\mathbb{T})$, say with constant $K_r$. (We will discuss verification of this hypothesis below.\footnote{The hypothesis \eqref{eq:unif} can be significantly weakened, see \cite{CarbCZ} and \cite{Seeger}.}) We shall then be able to conclude that 
\begin{equation}\label{eq:multthm}
\|\sum_{j\geq0} T_j f\|_p \leq C_p \|f\|_p
\end{equation}
for all $1 < p < \infty$, where $C_p = A_p^{-1} B_p K_{(\max\{p, p'\}/2)'}^{1/2}$.
Note that orthogonality already handles the case $p=2$.

To see this, we first note that inequality \eqref{eq:multthm} is essentially self-adjoint, meaning that if we can prove it for $p > 2$, then we can also prove it for $p < 2$. So assume $p >2.$ Then, by \eqref{eq:LP}, using the facts that $Q_k T_j = 0$ unless $j=k$ and that $Q_k$ and $T_k$ commute,
\[
\|\sum_{j\geq0} T_j f\|_p \leq A_p^{-1} \| \big( \sum_{k \geq 0}  |Q_k (\sum_{j \geq 0} T_j f)|^2\big)^{1/2}\|_p = A_p^{-1} \| \big( \sum_{k \geq 0}  |T_k Q_k f|^2\big)^{1/2}\|_p.
\]
Now we can calculate the expression 
$\| \big( \sum_{k \geq 0}  |T_k Q_k f|^2\big)^{1/2}\|_p$ using H\"older's inequality and its converse by
\[ \| \big( \sum_{k \geq 0}  |T_k Q_k f|^2\big)^{1/2}\|_p^2 = \sup_{\|w\|_r \leq 1} \int  \sum_{k \geq 0}  |T_k Q_k f|^2 w\]
where $ r = (p/2)'$. Fixing $w$ with $\|w\|_r \leq 1$, we have, by \eqref{eq:unif},
\[
\int  \sum_{k \geq 0}  |T_k Q_k f|^2 w \leq  \int \sum_{k \geq 0} |Q_k f|^2 \mathcal{M} w  \leq \| \big( \sum_{k \geq 0}  |Q_k f|^2\big)^{1/2}\|_p^2 \|\mathcal{M}w\|_r \leq K_r  \| \big( \sum_{k \geq 0}  |Q_k f|^2\big)^{1/2}\|_p^2,
\]
and we conclude, using \eqref{eq:LP} once more, that 
\[  \| \big( \sum_{k \geq 0}  |T_k Q_k f|^2\big)^{1/2}\|_p \leq K_r^{1/2} 
 \| \big( \sum_{k \geq 0}  |Q_k f|^2\big)^{1/2}\|_p \leq K_r^{1/2} B_p \|f\|_p,
 \]
which establishes the desired conclusion under the hypothesis \eqref{eq:unif}. 

\begin{remark}
(i) A very similar argument shows that under the same hypothesis \eqref{eq:unif} we also have $\|\big(\sum_{j\geq 0} |T_jf|^2\big)^{1/2}\|_p \leq C_p' \|f\|_p$. 

(ii) It is perhaps more traditional to see the above argument presented as a pointwise inequality of the form \[{\mathcal{Q}_1}(Tf)(x) \lesssim \mathcal{Q}_2(f)(x),\]
where $T = \sum_{j\geq 0} T_j$, and $\mathcal{Q}_1$ and $\mathcal{Q}_2$ are certain variants of $\mathcal{Q}$, enjoying the same $L^p$ inequalities as $\mathcal{Q}$ when $p \geq 2$, see \cite{Stein}. (Here and below we use the notation $\lesssim$ to denote the existence of a constant bounding the ratio of the left hand side to the right hand side.)

(iii) This argument crystallises the idea that the Littlewood--Paley theory allows us to derive conclusions about a {\em sum} of operators $\sum_j T_j$ from detailed uniform knowledge of the $T_j$ {\em individually}, as in \eqref{eq:unif}. 
\end{remark}

The uniform weighted assumption \eqref{eq:unif} is easily verified, for example, under the hypotheses\footnote{These are essentially quantitative uniform smoothness hypotheses on the Fourier multipliers $m_j$.} of the H\"ormander--Mikhlin multiplier theorem. In that case we simply take $\mathcal{M}w$ to be (a constant multiple of) the {\em Hardy--Littlewood maximal function} of $w$.
We turn next to the circle of ideas around maximal functions.

In \cite{LP1}, Littlewood and Paley were strongly motivated by potential applications to Lusin's problem concerning pointwise convergence of Fourier series.
While this was settled some thirty-five years later by Carleson \cite{Carleson} and Hunt \cite{MR238019}, in the early 1930's the question was still very much alive. Already it was understood that the correct way to control pointwise convergence behaviour in a quantitative way is with the use of appropriate {\em maximal functions}, of which the newly-born Hardy--Littlewood maximal function \cite{MR1555303} was and is the example {\em par excellence}. The (one-dimensional) Hardy--Littlewood maximal operator 
\[f \mapsto Mf(x) :=\sup_{r>0} \frac{1}{2r} \int_{B(x,r)} |f|\] 
is bounded on $L^p$ for $p>1$, and is of so-called weak-type $(1,1)$. This gives a quantitative version of Lebesgue's differentiation theorem stating that for $f$ in $L^1(\mathbb{T})$, the averages $ {(2r)^{-1}}\int_{B(x,r)} f $ converge to $f(x)$ for almost every $x$. Back in the context of pointwise convergence of Fourier series, the aim at the time was to reduce to a bare minimum any additional smoothness hypotheses to be placed on $f \in L^p$ for pointwise almost-everywhere convergence to hold, or, alternatively, to prove almost-everywhere convergence for all $f \in L^p$ for some especially interesting subsequences of partial sums of the Fourier series. Littlewood and Paley in the paper under discussion made significant contributions under both headings to what was then known.

For a function $f \in L^1(\mathbb{T})$, the $N$'th partial sum of its Fourier series is
\[ S_N f(x) = \sum_{|j| \leq N} \widehat{f}(j) e^{2 \pi ijx}.\]
The maximal function which controls the pointwise behaviour of $S_Nf(x)$ is given by
\[
S_\ast f(x) = \sup_{N \geq 1} |S_N f(x)|;\] 
the maximal function which controls the pointwise behaviour of the lacunary partial sums $S_{2^k}f(x)$ for $k \geq 0$ is given instead by
\[
S_\ast f(x) = \sup_{k \geq 0} |S_{2^k}f(x)|.\] 

One seeks to prove that (either version of) $S_\ast f$ is finite almost everywhere for all $f$ in as large a subclass $\mathcal{L}$ of $L^1$ as possible -- then pointwise almost-everywhere convergence for the same class will hold, via the fact that it always holds for sufficiently nice functions $f$ by classical results, and then a density argument. Kolmogorov \cite{zbMATH02586338} had shown that we cannot take $\mathcal{L}$ to be all of $L^1$ since there exist functions in $L^1$ whose Fourier series diverge almost everywhere, and indeed everywhere, but the question for $\mathcal{L}=L^p$ for $p >1$ was still open. And the best way to obtain almost-everywhere finiteness of $S_\ast f$ for $f \in \mathcal{L} \subseteq L^p$ is to show that we have an {\em a priori} inequality of the form
\[ \|S_\ast f\|_p \lesssim \|f\|_{\mathcal{L}}.\]

In this paper \cite{LP1} Littlewood and Paley announce:
\begin{theorem}\label{thm:conv} Suppose $1 < p \leq 2$ and $f \in L^p$.
\begin{enumerate}
\item [(i)] Then we have
\[ \| \sup_{k \geq 0} |S_{2^k} f|\|_p \lesssim \|f\|_p.\]
\item [(ii)] If in addition we have that $\{\widehat{f}(j) \log (|j| +2)^{1/p}\}_{j \in \mathbb{Z}}$ are the Fourier coefficients of some function $g \in L^p$, then 
\[ \|  \sup_{N \geq 1} |S_N f(x)|\|_p \lesssim \|f\|_p + \|g\|_p.\]
\end{enumerate}
\end{theorem}
We need to say a few words about the condition on $f$ in part (ii). It is to be thought of as a rather mild smoothness condition on $f$. This is because if $f$ and $f'$ are continuous, 
then $\widehat{f'}(k) = - 2 \pi i k  \widehat{f}(k)$, as an elementary calculation shows. Similarly, for $m \geq 1$, $\widehat{f^{(m)}}(k) = (- 2 \pi i k)^m  \widehat{f}(k)$. Moreover the $L^p$ norms of the functions with Fourier coefficients $(- 2 \pi i k)^m  \widehat{f}(k) \mbox{ and } ( 2 \pi  |k|)^m  \widehat{f}(k)$
are equivalent, because they are related by a Fourier multiplier of H\"ormander--Mikhlin type as described after Theorem~\ref{thm:Marc_primordial}. Defining the fractional derivative $D^sf$ of $f$ (where $s$ is now any non-negative real number) as the function (if it exists!) whose Fourier coefficients are $ (2 \pi  |k|)^s  \widehat{f}(k)$, it is thus natural to think of $\widehat{f}(j) \log (|j| +2)^{1/p}$ as representing the Fourier coefficients of a logarithmic variant (near $s =0$) of the fractional derivative $D^sf$.\footnote{Indeed, the Sobolev space $L^p_s$ can be defined to consist of those $f \in L^p$ such that $D^sf \in L^p$.}

Part (i) of this theorem in the case $p=2$ had already been observed by Hardy and Littlewood. It also holds for $2 < p < \infty$, but the pointwise convergence consequence in that case already follows from the case $p=2$. The pointwise convergence consequence of part (ii) in the case $p=2$ had been observed by Kolmogorov, Seliverstov and Plessner, but the maximal function bound was new both in this case and the case $1 < p <2$. The bound for the maximal function in part (ii) when $p>2$ is also given in the paper, but under a stronger than desired hypothesis on the Fourier coefficients, and Littlewood and Paley were clearly dissatisfied with it even at the time. Nevertheless, this stronger hypothesis was  weakened (and indeed removed altogether) only with Carleson's work in the 1960's.  

The proof of part (i) of this theorem is rather straightforward, but it has been enormously influential. The idea, which goes back to Kaczmarz and Zygmund,\footnote{Stein \cite{SteinZygmund} attributes the first appearance of quadratic expressions in harmonic analysis, albeit in the auxiliary role we are now discussing, to Kaczmarz and Zygmund in 1926, though this was perhaps not known at the time by Hardy, Littlewood and Paley.} is as follows: we can decompose $S_{2^k}$ into contributions $P_k$, which is an averaging operator at scale $2^{-k}$, corresponding to frequencies of size at most $2^k$, and $T_k$ for frequencies of size about $2^k$, in such a way that 
\[S_{2^k} = P_k + T_k,\]
where we have an easy pointwise control $ \sup_k |P_k f(x)| \lesssim Mf(x)$, and where $\{T_k\}$ satisfies the considerations of Remark (i) after Theorem~\ref{thm:Marc_primordial}. The key point is then that 
\[ \sup_k |S_{2^k} f| \leq  \sup_k |P_k f| + \sup_k |T_kf| \lesssim Mf + \big(\sum_k |T_kf|^2\big)^{1/2},\]
and the second inequality here is not foolish or wasteful precisely because of the disjoint supports of the Fourier multipliers corresponding to $\{T_k\}$.

This argument serves as the prototype of a vast array of applications of Littlewood--Paley theory to the analysis of maximal functions. The thrust of the method is that a quadratic expression mediates between -- or acts as a Tauberian condition linking -- the more difficult maximal function under study, and a more basic maximal function which is already understood. This reduces matters to establishing $L^p$ boundedness of the mediating quadratic expression. Hitherto, the study of maximal averaging operators in higher dimensions had relied upon establishing delicate geometric covering lemmas somewhat analogous to the well-known Vitali or Wiener covering lemmas for balls. But in many cases of interest, such as spherical averages (in which averages over balls appearing in the higher dimensional Hardy--Littlewood maximal operator are replaced by averages over spheres, see \cite{SteinBig}), it was not known how to do this. Mediation by quadratic expressions offered an entirely different perspective on such problems, even in the case $p=2$, and in the 1970's and 1980's led to spectacular success not only for the spherical maximal operator but also for (amongst many others) the study of differentiation in lacunary directions \cite{NSW}. The more basic, ``already-understood" maximal functions in these two particular cases are the classical Hardy--Littlewood maximal operator and the Hardy--Littlewood maximal operator with parabolic dilations (see below) respectively. 

While Theorem~\ref{thm:conv} has been superseded in strict terms by the seminal work of Carleson and Hunt \cite{Carleson}, \cite{MR238019}, the fundamental ideas behind it continue to remain highly influential to this day.

The paper under discussion contains no proofs, and no details of the application of Theorem~\ref{thm:LP} to Theorem~\ref{thm:conv} part (ii). These are given in the subsequent papers \cite{LP2} and \cite{LP3} which appeared after Paley's tragic death in a skiing accident in 1933. The paper does sketch, in broad strokes, the main lines of argument, which are inextricably rooted in the function theory of one complex variable, and which rely on various considerations which are available only in that setting. It makes implicit contact with other quadratic expressions such as the area function introduced by Lusin in 1930, and later explored by Marcinkiewicz and Zygmund. 




\section{What came next?}
While Littlewood--Paley theory was slow to enter the mainstream canon of harmonic analysis, its scope and influence have grown exponentially over the last nearly 100 years. Indeed, Zygmund regarded Littlewood--Paley theory as well ahead of its time (see the remarks in Bob Fefferman's Preface to the 3rd edition of Trigonometric Series \cite{Zygmund3}), and Zygmund, together with Marcinkiewicz, were among the first to recognise its extraordinary flexibility and potential. Indeed, Marcinkiewicz soon proved the Fourier multiplier theorem \cite{Marcinkiewicz1939} which bears his name, and which we briefly discussed above. With remarkable prescience he already formulated his theorem in the multi-dimensional and multi-parameter setting which we discuss below.

But beyond a few sporadic results, it was not really until the 1950's that Littlewood--Paley theory began to emerge within harmonic analysis as a general set of interlinked ideas with an inbuilt flexibility and resilience. (The first book-level treatment was given in the 1959 edition of Zygmund's classic {\em Trigonometric Series} \cite{Zygmund2}.) And to realise its potential, it was necessary to free the theory from its origins in complex variables, and give it a new identity. The first step on this path was to replace holomorphic functions by harmonic functions and Poisson integrals: this allowed, in principle, a development of the theory in higher dimensions. Littlewood--Paley theory was finally freed from all vestiges of complex function theory with the introduction of purely real-variable methods.  This came about as a result of the hugely significant and successful theory of Calder\'on and Zygmund from the early 1950's which addressed Singular Integral Operators, the primordial example of which is the Hilbert transform\footnote{The ubiquity of the Hilbert transform is partly explained by the fact that it takes a function $f$ on the real line to the function whose harmonic extension to the upper half plane is precisely the conjugate function of the harmonic extension of $f$.}  
\begin{equation}\label{eq:hilb}
Hf(x) = \frac{1}{\pi} \int_{\mathbb{R}} \frac{f(x-y)}{y}{\rm d}y
\end{equation} 
acting on functions in $L^p(\mathbb{R})$. It is beyond the remit of this article to discuss Calder\'on--Zygmund theory
in any detail, but suffice it to say that it was quickly recognised, primarily by Stein, that the Littlewood--Paley theory fell directly under its scope, and this real-variables perspective was fully developed in Stein's {\em Singular Integrals and Differentiability Properties of Functions} \cite{Stein} in 1970. The principal point here is that the manifestly non-linear mapping
\[ f \mapsto \mathcal{Q}f (x)= \big(\sum_j |Q_j f(x)|^2\big)^{1/2}
\]
can be regarded as the {\em linear} mapping 
\[ f \mapsto \{Q_j f\}_j\]
evaluated in the Hilbert space $\ell^2$, and, as such, can be treated directly\footnote{Strictly speaking one needs to work with mollified versions of the dyadic block operators $Q_j$ as below here, but the passage between the two versions is not difficult, especially within a vector-valued framework.}  within the vector-valued incarnation of Calder\'on--Zygmund theory.\footnote{Alternatively, one my make bounds on the {\em linear} operator $f \mapsto \sum_j \pm Q_j f$ uniformly in all choices of $\pm$, and then use Khintchine's inequality to pass to the quadratic expression.} At a stroke, Littlewood--Paley theory simultaneously became available in any setting in which Calder\'on--Zygmund theory made sense, in particular in euclidean spaces $\mathbb{R}^n$, and could now be based entirely on real-variable methods. This gave an enormous boost to its potential influence and range of applications. 

Euclidean spaces $\mathbb{R}^n$ have a much richer set of symmetries than the unit circle $\mathbb{T}$, and it quickly became apparent that the quadratic expressions of Littlewood--Paley theory should be formulated with these symmetries in mind. Thus, on $\mathbb{R}^n$ we have the dilation symmetries $x \mapsto tx$ for $t >0$, and to them we may associate the discrete and continuous quadratic expressions or ``square functions"
\[ (\sum_{j\in \mathbb{Z}} |Q_{2^j} f|^2)^{1/2} \mbox{ and }
\left(\int_0^\infty |Q_t f|^2 \frac{{\rm d} t}{t}\right)^{1/2}
\]
as follows. We begin with a smooth bump function $\psi$, with mean zero, and which is chosen to satisfy
\begin{equation}\label{eq:smoothid}
\sum_{j\in \mathbb{Z}} |\widehat{\psi}(2^j \xi)|^2 \equiv 1 \mbox{ 
or } \int_{0}^\infty |\widehat{\psi}(t\xi)|^2 \frac{{\rm d} t}{t}\equiv 1 
\end{equation} 
where $\widehat{\cdot}$ denotes the Fourier transform on $\mathbb{R}^n$ defined by 
\[ \widehat{f}(\xi) = \int_{\mathbb{R}^n} {f}(x) e^{-2 \pi i x \cdot \xi} {\rm d} x.\]
The Fourier transform is an isometry of $L^2$ according to Plancherel's theorem. Construction of such $\psi$ can be achieved by means of simple partition of unity arguments. We now let 
\[ 
Q_t f = \psi_t \ast f = \frac{1}{t^{n}} \psi \left(\frac{\cdot}{t}\right) \ast f.
\]
where $\ast$ denotes convolution (which is converted to multiplication by the Fourier transform). Plancherel's theorem together with \eqref{eq:smoothid} immediately gives the Pythagorean identities
\[ \|  (\sum_{j\in \mathbb{Z}} |Q_{2^j} f|^2)^{1/2}\|_2 = \|f\|_2
\mbox{ and } \|\left(\int_0^\infty |Q_t f|^2 \frac{{\rm d} t}{t}\right)^{1/2}\|_2 = \|f\|_2\]
respectively.\footnote{The reader will note that we have switched notation from $Q_j$ earlier to $Q_{2^j}$ here. This is to explicitly reflect the dilation {\em group} structure, and for consistency between the discrete and continuous formulations. We revert to the original notation in Section~\ref{sec:exotic}.} 

One then has, as a direct consequence of Calder\'on--Zygmund theory, the fundamental relations
\[ \|  (\sum_{j\in \mathbb{Z}} |Q_{2^j} f|^2)^{1/2}\|_p \approx  \|f\|_p \mbox{ and } \|\left(\int_0^\infty |Q_t f|^2 \frac{{\rm d} t}{t}\right)^{1/2}\|_p \approx \|f\|_p\]
for $1 < p < \infty$. 

\begin{remark}
(i) The discrete and continuous versions are to be regarded as essentially interchangeable twin manifestations of the same phenomenon. 

(ii) In the presence of the Pythagorean identities the reverse inequalities follow from the direct ones via polarisation.

(iii) At this level of generality (and indeed lower levels), relevant quadratic expressions come in many flavours, (such as $\mathcal{Q}_1$ and $\mathcal{Q}_2$ mentioned in Section~\ref{sec:secondsection}), but we have suppressed discussion of these variants in the interests of clarity and brevity. 

(iv) Imaginative constructions of such quadratic expressions extend the application of the theory to operators whose oscillatory nature lies well outside the scope of classical Calder\'on--Zygmund theory, see for example \cite{BB}.
\end{remark}

In addition to the standard one-parameter dilations, there are many other families of symmetries of euclidean space. One such is the one-parameter {\em non-isotropic} dilations $x \mapsto (t^{a_1} x_1, \dots, t^{a_n} x_n)$ for $t > 0$ and fixed exponents $a_1, \dots, a_n >0$. While these can be treated by variants of the Calder\'on--Zygmund theory adapted to such dilations, the case of more general {\em multi-parameter} dilations $x \mapsto (t_1 x_1, \dots , t_n x_n)$ for {\em independent} parameters $t_1, \dots , t_n >0$ has a somewhat different character; but some results can be relatively easily established essentially via iterating the one-dimensional results. This is the framework within which the original Marcinkiewicz multiplier theorem fits. There are also dilational actions which have a hybrid nature, located between these two cases. All these have been thoroughly studied by Ricci and Stein and their collaborators, see for example \cite{RicciStein}.

On the other hand, systematic use of Littlewood--Paley (and vector-valued) theory in the non-isotropic setting (see especially \cite{MR837527}, \cite{Christ} and \cite{RicciStein}) feeds back into classical singular integral theory, and, rather remarkably, leads to a self-improving structural broadening of its aegis. For example, operators on $L^p(\mathbb{R}^n)$ such as the Hilbert transform along the curve $\Gamma(t) = (t^{a_1}, \dots, t^{a_n})$ given by 
\[ H_\Gamma f(x) = \frac{1}{\pi} \int_{\mathbb{R}} f(x-\Gamma(t))\frac{{\rm d}t}{t}
\] 
and its associated maximal function 
have a singularity which is far too severe to allow direct treatment by the classical non-isotropic Calder\'on--Zygmund theory. Nevertheless, they do fall under this expanded scope, and have been studied by Christ, Fabes, Nagel, Ricci, Rivi\`ere, Stein, Wainger and others, (see for example \cite{SW}, \cite{Christ}, \cite{RicciStein}), both in the euclidean and more generally the homogeneous group settings.

When one passes to non-dilational symmetries such as translations or rotations, however, the corresponding Littlewood--Paley quadratic expressions fall outside the direct purview of Calder\'on--Zygmund theory. They may nonetheless still contain a great deal of geometrical content relating to almost-orthogonality and be highly relevant for applications. See especially Section~\ref{sec:exotic} below.

An example of the flexibility afforded by the real-variables formalism of Littlewood--Paley theory is the ease with which it characterises whole scales of function spaces including Sobolev spaces, Bessel potential spaces, Triebel--Lizorkin spaces and Besov spaces. (Here we work with standard dilations of $\mathbb{R}^n$.) Indeed, the (homogeneous) Triebel--Lizorkin spaces $F^{p,q}_s$ -- which when $q=2$ correspond to the Bessel potential and Sobolev spaces $L^p_s$, and in particular, for $s=0$, to plain vanilla $L^p$ -- are characterised by finiteness of the expression
\[ \|(\sum_{j \in \mathbb{Z}} |2^{js}Q_{2^j}f|^q)^{1/q} \|_p.\] 
On the other hand, the (homogeneous) Besov spaces $B^{p,q}_s$ are obtained by taking the $L^p$ and $\ell^q$ norms in the ``wrong" order, and are characterised by finiteness of 
\[ (\sum_{j \in \mathbb{Z}} 2^{jsq}\|Q_{2^j}f\|_p^q)^{1/q}.\] 
Amusingly, another instance of taking norms in the the ``wrong" order is discussed below in Section~\ref{sec:exotic}.

With its real-variables foundation firmly in place, Littlewood--Paley theory was now ripe to be developed to its fullest extent.

\section{Quasi-abstract Littlewood--Paley theories}
What are the basic ingredients needed to formulate Littlewood--Paley functionals and develop Littlewood--Paley theory? Despite first appearances, they can be seen to have not so much to do with Fourier Analysis {\em per se}. Most importantly, what is essential is a {\em family of  averaging operators} $P_t$ on different ``scales" given by $t>0$, each of which enjoys further natural (or even necessary) properties, and which are related at different scales in some appropriate way. The operators $Q_{2^j}$ and $Q_t$ are obtained by taking respectively successive differences (i.e. $Q_{2^j} = P_{2^j} - P_{2^{j-1}}$ ) or logarithmic derivatives (i.e. $Q_t = t \frac{\partial P_t}{\partial t}$) of these averaging operators $P_{t}$.

The first level of abstraction takes place in the realm of conditional expectation and martingales (and thus martingale differences) on probability spaces or $\sigma$-finite measure spaces. This development was to a very large extent presaged by Paley's 1932 paper \cite{MR1576148}, and culminated in the 1960's and 1970's with seminal works of Burkholder, Davis and Gundy \cite{MR208647}, \cite{MR221573} and \cite{MR400380}. (Interestingly, Burkholder \cite{MR208647} credits Gundy as being the first to recognise that Paley's paper was susceptible to substantial generalisation.) The central result is the analogue of Theorem~\ref{thm:LP} in this setting, 
the equivalence of $L^p$ norms of a martingale with those of its martingale difference square function. This subsumed an earlier result of Marcinkiewicz and Zygmund \cite{MZ} with a similar flavour concerning independent random variables with mean zero. In conjunction with the remark after Theorem~\ref{thm:LP}, that result explains why lacunary Fourier series are sometimes said to behave in the same way as sums of independent random variables. The Marcinkiewicz--Zygmund result in turn subsumes Khintchine's inequality. It is the results of 
Burkholder, Davis and Gundy which unlock all of the potential of Littlewood--Paley theory from Fourier Analysis into the stochastic and indeed more general settings.

Simultaneously with Burkholder's work, Stein was beginning to develop Littlewood--Paley theory on compact Lie groups. He realised that the correct setting for this was that of symmetric diffusion semigroups (mimicking the Poisson semigroup in euclidean spaces). In its most general formulation, his theory relies on the Hopf--Dunford--Schwartz maximal ergodic theorem, together with estimates for quadratic expressions which are (nontrivially!) obtained from the corresponding ones in the setting of martingale differences. This, together with many other developments concerning the essentially holomorphic nature of the semigroup setting and its application to holomorphic spectral multipliers, is set out in \cite{SteinTHA}. Later, Cowling \cite{Cowling1}, \cite{Cowling2} offered an alternative and perhaps simpler perspective -- based upon transference -- on this circle of ideas.

These developments -- in both the martingale and the semigroup settings -- have had a profound influence which is worthy of far more extensive discussion than we can give here. For now we wish merely to note that the power that abstraction brings can be exploited and re-channelled to give fresh insight into even the most classical of settings, with spectacular results, particularly when it comes to the quantitative study of bounds. Indeed, the semigroup approach leads to dimension-free bounds for the Hardy--Littlewood maximal operator defined over balls, and even over arbitrary convex bodies in $\mathbb{R}^n$ (see \cite{SteinBig}, \cite{Bourgain}, \cite{CarbConvex}), and martingale methods lead to the best constant in the $\ell^p(\mathbb{Z})$ bounds for the discrete version of the Hilbert transform \eqref{eq:hilb}. 
In this latter context see especially the survey article of Ba\~nuelos \cite{ban}.

\section{Vector-valued analysis and weighted inequalities}
Vector-valued analysis arising out of Littlewood--Paley theory brings a new flexibility of perspective. In short, instead of considering scalar-valued operators $T$, we consider operators taking possibly vector-valued inputs to possibly vector-valued outputs. The whole machinery of singular integral operators adapts readily to this situation. As we have seen above, Littlewood--Paley theory neatly fits this paradigm, (in which case the domain consists of scalar-valued functions and the codomain of Hilbert space-valued functions). The Marcinkiewicz--Zygmund result saying that any bounded linear operator on $L^p$ automatically extends to a bounded linear operator on $L^p(\ell^2)$ with the same norm is an early example of the vector-valued philosophy. Apart from anything else, vector-valued inequalities now form an essential part of the working toolkit of modern Fourier analysis. There is nothing limiting us to considering only Hilbert spaces for our vectors, and indeed when considering general Banach spaces, properties such as UMD (unconditional martingale differences) turn out to play a key role (for example, the Hilbert transform is bounded on $L^p$ with values in a Banach space $X$ if and only if $X$ is UMD, see \cite{Bourg1}, and Littlewood--Paley inequalities are valid in $L^P(X)$ spaces if and only if $X$ is UMD, see  \cite{MR752501}).

Keeping things down to earth, however, suppose that we have a family of operators $T_j$ satisfying
\begin{equation}\label{eq:unif_2}
\sup_j \int |T_jf|^2 w \lesssim \int |f|^2 \mathcal{M}w
\end{equation}
where $w \mapsto \mathcal{M}w$ is some (possibly non-linear) non-negative operator which is bounded on $L^r$. (We saw this above as \eqref{eq:unif}.) Then a very simple argument using H\"older's inequality shows us that we have
\begin{equation}\label{eq:vector}
\| (\sum_j  |T_jf_j|^2)^{1/2}\|_p \lesssim \| (\sum_j |f_j|^2)^{1/2}\|_p
\end{equation}
for $p>2$ such that $(p/2)' = r$. It is a quite remarkable fact that the converse to this holds: under rather general conditions, an inequality of the sort \eqref{eq:vector} necessarily implies that some inequality of the sort \eqref{eq:unif_2} holds. This principle -- that vector-valued inequalities are equivalent to weighted inequalities -- is quite general. But the catch is that the argument is non-constructive: even if we know perfectly what are the $T_j$, we in general have no idea what will be the mapping $w \mapsto \mathcal{M}w$. Even at the level of a single $T_j$ this remains remarkable, and can be captured in the maxim that ``there is no such thing as $L^p$ -- only $L^2$ with weights".

The book \cite{GCRdeF} by Garc\'{i}a-Cuerva and Rubio de Francia contains an extensive introduction to these topics and their interplay with Littlewood--Paley theory, covering up to the mid-1980's. The theory has continued to develop apace since then, not least in the context of non-commutative and operator-valued harmonic analysis, building on the martingale perspective. See \cite{LuPiqPis} and \cite{PisierXu}, and, for some more recent developments, \cite{KLW}.

\section{Wavelets and the Calder\'on reproducing formula} 
Thanks largely due to the incisive insights of Yves Meyer, for which he was awarded the Abel Prize, there has arisen since the 1980's a vast literature on wavelets and their relatives, to which we cannot possibly do justice. Three excellent references are \cite{Meyer}, \cite{Daubechies} and \cite{HernandezWeiss}. Beyond this, the applications of wavelets across science and engineering are countless.

The mathematical origins of wavelets lie in the Calder\'on reproducing formula
\begin{equation*}\label{eq:Calderon}
f(x)= 
\sum_{j\in \mathbb{Z}} {\psi}_{2^j} \ast \tilde{\psi}_{2^j} \ast f (x) = \sum_{j \in \mathbb{Z}} \int_{\mathbb{R}^n}\psi_{2^j}(x-y)  \tilde{\psi}_{2^j} \ast f(y) {\rm d}y
= \sum_{j \in \mathbb{Z}} \int_{\mathbb{R}^n}\psi_{2^j}(x-y)  \langle f, \tilde{\psi}_{2^j}(\cdot - y)\rangle {\rm d} y
\end{equation*}
which dates from the early 1960's and is essentially a re-statement of identity \eqref{eq:smoothid} (where $\tilde{\psi}(x) = \overline{\psi}(-x)$). Assume that $\psi$ is even and real-valued. Since $\psi_{2^j}$ has scale $2^j$, we may morally replace the $y$-integral here by a sum over $2^j$-spaced lattice points, and it becomes 
\[\sum_{j \in \mathbb{Z}} \sum_{k \in \mathbb{Z}^n} 2^{jn}\psi_{2^j}(x-2^j k)  \langle f, {\psi}_{2^j}(\cdot - 2^j k)\rangle = \sum_{j \in \mathbb{Z}, k \in \mathbb{Z}^n} \langle f, \Psi_{j, k} \rangle \Psi_{j,k}(x)
\]
where $\Psi_{j,k}(x) := 2^{-jn/2} \psi(2^{-j}x- k)$ is $L^2$-normalised. Here, $f \mapsto \{\langle f, \Psi_{j,k}\rangle\}_{j,k}$ is the wavelet transform of $f$, taking $f$ to its wavelet coefficients, and the formula represents the wavelet decomposition of $f$ into sums of wavelets $\{\Psi_{j,k}\}_{j, k}$. Realisation of this decomposition is related to sampling, and, for carefully chosen $\psi$, the representation of $f$ can be made rigorous. 

This is all naturally adapted to $L^2$, and the the $\{\Psi_{j,k}\}_{j, k}$ may be chosen to form a Riesz (or in some cases even an orthonormal) basis for $L^2$, so that $\|f\|_2 \approx (\sum_{j \in \mathbb{Z}, k \in \mathbb{Z}^n} |\langle f, \Psi_{j, k} \rangle|^2)^{1/2}$. 

However, the huge advantage of this approach is that it is also well-adapted to $L^p$ and many other natural function spaces. Indeed, the the $\{\Psi_{j,k}\}_{j, k}$ may be chosen to form an unconditional basis of $L^p$ for $1 < p < \infty$, and furthermore to satisfy the Littlewood--Paley-type relations
\begin{equation}\label{eq:LPwavelet}
\|f\|_p 
\sim \| (\sum_{j \in \mathbb{Z}, k \in \mathbb{Z}^n} |\langle f, \Psi_{j, k} \rangle|^2 |\Psi_{j,k}|^2)^{1/2}\|_p 
\sim \| (\sum_{j \in \mathbb{Z}, k \in \mathbb{Z}^n} 2^{-nj}|\langle f, \Psi_{j, k} \rangle|^2 \chi_{D_{jk}})^{1/2}\|_p
\end{equation}
where $D_{jk}$ is the dyadic cube of sidelength $2^j$ anchored at $2^{j}k$. Thus, in particular, all the information about the size of an $f \in L^p(\mathbb{R}^n)$ can be read off from just the {\em size} of its wavelet coefficients: this is in stark contrast to what happens in classical Fourier analysis, where the size of the Fourier coefficients is far from the full story; moreover it makes computation in the wavelet setting much more stable than in its Fourier counterpart.

Singular integral operators of both convolution and non-convolution form are well-quasi-diagonalised in $L^p$ by wavelet expansions. In fact the proof of the Littlewood--Paley relations \eqref{eq:LPwavelet} proceeds, unsurprisingly, by means of vector-valued Calder\'on--Zygmund theory.

A related development is the study of paraproducts which arise when the Calder\'on reproducing formula is used to write a product $fg$ as 
\[ fg = (\sum_i \Delta_{2^i} f)(\sum_j \Delta_{2^j} g) = \sum_j  \Delta_{2^j} f \Delta_{2^j} g + \sum_{i<j}  \Delta_{2^i} f \Delta_{2^j} g  + \sum_{j<i}  \Delta_{2^i} f \Delta_{2^j} g,
\]
where $\Delta_{2^j}$ is convolution with $Q_{2^j} \ast \tilde{Q}_{2^{j}}$. Setting $P_{2^j} = \sum_{i: i <j} \Delta_i$ this becomes
\[ fg = \sum_j  \Delta_{2^j} f \Delta_{2^j} g + \sum_{j}  P_{2^j} f \Delta_{2^j} g  + \sum_{j}  \Delta_{2^j} f P_{2^j} g.
\]
The second and third terms arising on the right-hand side here are examples of paraproducts, so-called because they mimic classical products but with a limited interaction between Fourier frequencies. 

One of many natural that questions one may pose about such decompositions is whether there are estimates for each of the three terms which allow us to recover the H\"older inequality $\|fg\|_p \lesssim \|f\|_q \|g\|_r$ where $1/p = 1/q + 1/r$.
The first term does not present difficulties, and can be treated by Cauchy--Schwarz and then Littlewood--Paley theory directly. The paraproduct terms can be treated using Carleson measures and fall under the general theory of non-convolution singular integrals of David, Journ\'e and Semmes going by the name of $T1$ and $Tb$ theorems. See \cite{SteinBig}, \cite{DJ}, \cite{DJS}, \cite{Hofmann}. Notably, a general singular integral operator of the David--Journ\'e--Semmes class can be written as a sum of three terms, one of which has special cancellation properties such as $T1 = T^\ast1 =0$, and the other two of which are paraproducts.

Paraproducts and their relatives also arise in non-linear PDE (in Bony's paradifferential calculus) and are used to establish versions of the Leibniz formula (i.e. the product rule) for fractional derivatives. 


\section{Almost-orthogonality principles beyond Littlewood--Paley}
Littlewood--Paley theory is, in essence, a theory of almost-orthogonality for functions in and operators on $L^p$ when $p\neq 2$. However, there are operators on $L^2$ for which one may struggle to do a spectral analysis, and for which an almost-orthogonality theory based on real variables, and 
going beyond that of Littlewood and Paley, would be very useful. Such a principle is given by the Cotlar--Stein lemma:

\begin{theorem}[Cotlar--Stein, see \cite{SteinBig}] Suppose that $H$ is a Hilbert space and that $T_j \in \mathcal{L}(H)$ satisfy
\[ \|T_j T_k^\ast\| + \|T_j^\ast T_k\| \lesssim 2^{- \epsilon |j-k|}\]
for all $j, k \in \mathbb{Z}$. Then $T := \sum_j T_j \in \mathcal{L}(H)$. 
\end{theorem}
The proof is by a clever but elementary combinatorial argument. The hypothesis is that we have $\sup_j \|T_j\| < \infty$, plus some almost-orthogonality between the different $T_j$. Generally one has an operator $T$ at hand whose boundedness one wants to establish; the art is in finding a natural  decomposition $T = \sum_j T_j$ (if one exists) to which one may apply the Cotlar--Stein lemma. In the case of the Hilbert transform \eqref{eq:hilb}, for example, it is natural to take $T_jf(x) = \frac{1}{\pi} \int_{|y| \sim 2^j} \frac{f(x-y)}{y}{\rm d}y$ (of course we do not need to do this, because we already know the Hilbert transform is bounded). Even when a spectral decomposition of $T$ is not readily available, such real-variables decompositions often work provided that $T$ satisfies some cancellation and smoothness properties. Boundedness on $L^2$ of various classes of pseudodifferential operators of Calder\'on--Vaillancourt-type was amongst the first successes of the Cotlar--Stein lemma, see \cite{SteinBig}. 

The hypothesis can be verified in certain concrete situations when $H = L^2$ as follows. By Schur's lemma it suffices to show that, with $\|T\|_{p-p}$ denoting the $L^p$ operator norm, 
\[ \|T_j T_k^\ast\|_{1-1} + \|T_j^\ast T_k\|_{1-1} \lesssim 2^{- \epsilon |j-k|}\]
for all $j, k \in \mathbb{Z}$, and this may be directly verified if $T$ has appropriate cancellation and its integral kernel has some smoothness. This works in the case of classical non-convolution singular integral operators $T$, when, say, $T1=T^\ast 1 =0$. But even when the integral kernel of $T$ is only a measure, and so has no classical smoothness, the technique can still be made to work. Let $S=T_jT_k^\ast$ or $T_j^\ast T_k$. The idea is that instead of estimating $\|S\|$ directly by $\|S\|_{1-1}$, we instead estimate $\|S\|^2 =\|S^\ast S\|$ by $\|S^\ast S\|_{1-1}$ (or by $\|SS^\ast\|_{1-1}$). This has more chance of working because the effect of the composition of $S^\ast$ with $S$ is typically to spread out and smooth out the integral kernel, as in the case of repeated convolution. The more repetitions of this strategy one carries out, the more likely one is to succeed! For an early manifestation of this idea see the work of Christ, \cite{Christ}. 

Following a suggestion of Michael Cowling, the Cotlar--Stein lemma can be used in conjunction with Littlewood--Paley theory to give further $L^p$-almost-orthogonality principles which are applicable in non-commutative and non-translation-invariant settings, see for example \cite{CarbCotlar}. Moreover, as observed by Christ, such principles may be further adapted to give almost-orthogonality principles for maximal functions, the roots of which lie in the work of Nagel, Stein and Wainger \cite{NSW} on differentiation in lacunary directions. See \cite{MR837527} and \cite{CarbConvex}; see also \cite{MR1809116}.

\section{Littlewood--Paley theory in PDE and Geometric Measure Theory}\label{sec:GMT}

Littlewood--Paley theory and quadratic expressions -- especially those such as Lusin's area function -- were, in their nascence, intimately connected to complex function theory and to boundary value problems for harmonic functions on the upper half-plane. It it therefore unsurprising that they continue to play a key role in the more modern and refined study of the solvability of boundary value problems for elliptic and parabolic PDE, including those with rough coefficients, and on domains with rough boundaries where the boundary data is in $L^p$. The quadratic expressions which occur are, naturally, often more closely related to the Lusin area function than to the more Fourier-analytic square functions that we have focused on in this article. However, Littlewood--Paley expressions associated to the diffusion semigroup generated by the uniformly elliptic operator $-{\rm div} (A \nabla)$ are also important. The survey article of Dindo\v{s} and Pipher \cite{DindosPipher} is a good reference for the modern theory for uniformly elliptic PDE with rough coefficients in Lipschitz domains in which quadratic expressions feature strongly.

\medskip
The analogue of the Hilbert transform in the classical setting becomes a singular integral -- a Cauchy-type integral -- associated to a rough (hyper-)surface,\footnote{It is to be noted that this is a quite different generalisation from $H_\Gamma$ discussed above.} and the study of such involves the development of singular integral theory going beyond the treatment of David, Journ\'e and Semmes. Calder\'on established the boundedness of the Cauchy integral along curves with small Lipschitz constant, and Coifman, McIntosh and Meyer did so for arbitrary Lipschitz curves. David, Jones and Semmes pioneered the study of still rougher curves and surfaces. The geometric measure theory of such surfaces begins to play a critical role: boundedness of such a singular integral operator turns out to be closely related to uniform rectifiability of the underlying surface. Quadratic expressions, along with Carleson measures, feature throughout the development, including in the celebrated solution of the Kato square root problem by Auscher, Hofmann, Lacey, McIntosh and Tchamitchian, pushing the theory of singular integrals to its limits, see for example \cite{Hofmann}; and also Tolsa's solution of the Painlev\'e problem\footnote{The Painlev\'e problem was one of the initial motivations behind the study of the Cauchy integral on Lipschitz graphs and uniformly rectifiable sets.}, see \cite{Tolsa}. See also \cite{Verdera}, Sections 7 and 8.

\medskip
Quadratic expressions involving Peter Jones' $\beta$-numbers and their relatives are central in the study of uniform rectifiability. They also give an analytic perspective on structure in large data sets -- in particular measuring clustering of data around rectifiable sets. See \cite{DavidSemmes} (Section I.1.3) which directly addresses the analogies between Littlewood--Paley techniques in analysis and in geometry. There is also a very nice survey \cite{Mattila} by Mattila on rectifiability, see especially Chapter 5 on uniform rectifiability. 
See \cite{MR3952696} for connections between the PDE perspective and $\beta$-numbers in higher-codimensional settings.

There is now a vast literature on this area which we can make no attempt to summarise; partly for this reason this section is thus perhaps even more impressionistic than the others. Amongst the contributors in the current era are Azzam, D\k{a}browski, David, Dindo\v{s}, Feneuil, Hofmann, Kenig, Li, Martell, Mattila, Mayboroda, Mourgoglou, Pipher, Tolsa, Toro, Verdera, Villa and their collaborators.

\section{Exotic Littlewood--Paley decompositions}\label{sec:exotic}

\subsection{Equally-spaced intervals} 

Carleson first considered non-standard Littlewood--Paley expressions where the collection of dyadic blocks (which is invariant under dilations) is replaced by the collection of disjoint translates of a fixed interval. Thus, in the setting of $\mathbb{R}$, he considered 
\[ f \mapsto \mathcal{Q}f = (\sum_{j \in \mathbb{Z}} |Q_jf|^2)^{1/2}\]
where now (reverting to the notation around $j$ of earlier sections)
\[ \widehat{Q_j f}(\xi) = \chi_{[j-1/2, j+1/2)}(\xi) \widehat{f}(\xi)\]
or its smooth version
\[\widehat{Q_j f}(\xi) = \phi(\xi-j) \widehat{f}(\xi) \]
where $\phi$ is a smooth normalised bump function. 

It is easily seen that any inequality of the form $\|\mathcal{Q} f\|_p \lesssim \|f\|_p$ for $p < 2$ with $N$ operators $Q_j$ fails by a factor $\sim N^{\frac{1}{p} - \frac{1}{2}}$, (and hence that there is no inequality of the form $\|f\|_p \lesssim \|\mathcal{Q} f\|_p$ for $p > 2$), but Carleson showed that  $\|\mathcal{Q} f\|_p \lesssim \|f\|_p$  does indeed hold for $p\geq 2$. Later, C\'ordoba and separately Rubio de Francia gave alternative arguments for this, with Rubio de Francia's argument \cite{RdeFequal} being especially informative: using nothing more than Parseval's identity and basic properties of the Fourier transform, he showed that for the smooth version of $\mathcal{Q}$, there is essentially an identity (modulo negligible error terms)
\[ \mathcal{Q} f(x) \approx \left( \Phi \ast |f| ^2(x)\right)^{1/2}\]
where $\Phi$ is a suitably chosen normalised non-negative bump function. Anything that can be said about $\mathcal{Q}$ can be readily derived from this relation. Perhaps surprisingly, these ``equally-spaced" Littlewood-Paley inequalities turn out to be important in problems concerning $L^p$-Fourier inversion in $\mathbb{R}^n$, see below.

\subsection{Arbitrary intervals -- the Rubio de Francia Littlewood--Paley inequality}
Shortly thereafter, Rubio de Francia \cite{RdeFarb} further surprised the harmonic analysis community when he proved that one may replace the congruent intervals $[j-1/2, j+1/2)$ considered above by {\em any disjoint collection of intervals whatsoever}, and still conclude that  
$\|\mathcal{Q} f\|_p \lesssim \|f\|_p$ for $2 \leq p < \infty$. Moreover, in this case, his simple identity is replaced by a pointwise inequality
\[ (\mathcal{Q} f)^\#(x) \lesssim  \left( M|f| ^2(x)\right)^{1/2}\]
where $\#$ represents the Fefferman--Stein sharp function which appears in the definition of the space BMO and which satisfies $\| F\|_p \lesssim \| F^\#\|_p$ for $1 < p < \infty$. (See \cite{SteinBig}.) Remarkably, Rubio de Francia's pointwise inequality is proved by establishing kernel estimates for $f \mapsto \mathcal{Q}(f)$ as a (non-symmetric) vector-valued Calder\'on--Zygmund operator; this is in spite of the fact that an arbitrary collection of disjoint intervals contains no vestiges of dyadic block or dilation structure. 

Rubio de Francia's inequality was later extended to arbitrary families of disjoint axis-parallel rectangles in $\mathbb{R}^n$. See \cite{Journe}, \cite{Lacey} and also \cite{Roncaletal}. 

\subsection{Angular decompositions -- reverse inequalities and curvature}
It is natural also to consider Littlewood--Paley expressions corresponding to angular decompositions of $\mathbb{R}^n$. 

Sticking to $\mathbb{R}^2$ for the moment, amongst the more natural of these are decompositions into lacunary angular sectors $\Gamma_j = \{(\xi_1, \xi_2) \in \mathbb{R}^2 \, : \, 2^j \leq |\xi_2|/|\xi_1| < 2^{j+1}\}$ for $j \in \mathbb{Z}$, and into equiangular sectors $\Delta_j = \{(\xi_1, \xi_2) \in \mathbb{R}^2 \, : \, 2 \pi j/N \leq \arctan(|\xi_2|/|\xi_1|) < 2 \pi(j+1)/N\}$ for $0 \leq j < N$.

Notice that the parabolic dilation $\xi \mapsto (2^k\xi_1, 2^{2k} \xi_2)$ maps $\Gamma_j$ to $\Gamma_{j+k}$, and this suggests that the Littlewood--Paley theory for lacunary angular sectors might be largely accessible within the established theory. Indeed, the Littlewood--Paley inequalities for smoothed versions of the lacunary sectors $\Gamma_j$ are as expected, and are a straightforward consequence of the Marcinkiewicz multiplier theorem. If we want to consider quadratic expressions corresponding to the rough $\Gamma_j$, we will need some auxiliary vector-valued and maximal inequalities corresponding to differentiation in lacunary directions \cite{NSW}.

The theory for the equiangular sectors $\Delta_j$, or their smooth variants, where we see an underlying rotational action, is more challenging and interesting. One does not expect to have estimates of the form $\|\mathcal{Q} f\|_p \lesssim \|f\|_p$ for $p < 2$ for the same reasons that the inequalities for equally-spaced intervals fail. But when $p > 2$, ``pure" estimates with bounds independent of $N$ are not available either, due to Kakeya set examples.\footnote{More on this below.} However, there is a dimension-dependent window of $p$'s with $p > 2$ in which one might hope for bounds of the form $\|\mathcal{Q} f\|_p \lesssim \|f\|_p$ with at most subpolynomial or logarithmic failure in $N$. In two dimensions this window is $2 \leq p \leq 4$, and in higher dimensions it is $2 \leq p \leq 2n/(n-1)$. C\'ordoba \cite{Cordoba} established the positive result in the two-dimensional case, while the higher-dimensional case awaits resolution of the Kakeya maximal conjecture, amongst other things, and remains open.

Just when it seems we have said everything that can currently be said, matters take an interesting turn. One of the most powerful drivers for research in harmonic analysis over the last fifty years has been the question of inverting the Fourier transform for functions in $L^p(\mathbb{R}^n)$ (as opposed to the almost-everywhere questions of Sections~\ref{sec:intro} and \ref{sec:epilogue}). The fundamental question here is, with 
\begin{equation}\label{eq:FIRn}
S_Rf(x) = \int_{|\xi| \leq R} \widehat{f}(\xi) e^{2 \pi i x \cdot \xi} {\rm d} \xi,
\end{equation}
whether $S_Rf \to f$ in $L^p$-norm for all $f \in L^p(\mathbb{R}^n)$, (say for $p \leq 2$), which, with a little elementary functional analysis, is easily seen to be equivalent to boundedness on $L^p(\mathbb{R}^n)$ of the operator $S_R$, or indeed of $S_1$. Fefferman \cite{Feffermanball} astounded the world of harmonic analysis when he showed that this boundedness failed when $n \geq 2$, except when $p=2$. And moreover what is responsible for this failure is the rough singularity of the multiplier $\chi_{|\xi|\leq 1}$ for $S_1$ at $|\xi| =1$, together with the fact that the normal directions to the sphere fill an open set. This allowed him to construct counterexamples based on Kakeya sets. For this reason, and others, it is natural to decompose the rough disc multiplier $\chi_{\{|\xi| <1\}}$ into smooth (Littlewood--Paley!) dyadic pieces $\phi_{2^j}$ for $j \geq 0$, each supported in an annulus of width $\sim 2^{-j}$ and distant $\sim 2^{-j}$ from the unit sphere, and to ask about $L^p$ bounds for the corresponding operators $S^{2^{-j}}$.\footnote{The question of $L^p$ orthogonality {\em between} the different $S^{2^{-j}}$ is an interesting matter which we will pursue elsewhere.} Obviously we have $\|S^{2^{-j}}\|_{2-2} \lesssim 1$ uniformly in $j$. Fix $j$, and let $\delta = 2^{-j}$. It is also natural to break up each $\phi_\delta$ into $\sim \delta^{-(n-1)/2}$ equiangular smooth pieces $\phi_\nu$, each supported in a slab of size $\sim \delta \times \delta^{1/2} \times \dots \times \delta^{1/2}$ contained
in the $\delta$-neighbourhood of the sphere, each of which has corresponding operator $S_\nu$ which is easily seen to satisfy $\|S_\nu \|_{p-p} \lesssim 1$ for $1 \leq p \leq \infty$. One quickly concludes using interpolation and the triangle inequality that $\|S^{\delta}\|_{p-p} \lesssim \delta^{-(n-1)|\frac{1}{p} - \frac{1}{2}|}$. However, the use of the triangle inequality here turns out to be a very poor choice, and we can do much better. Indeed, there is a huge amount of orthogonality {\em between} the different $S_\nu$ (remember $\delta$ is fixed here) which we have not thus far exploited. This orthogonality finds its formulation, when $n=2$, in the reverse Littlewood--Paley inequality
\[ \|\sum_\nu S_\nu f\|_4 \lesssim  \|(\sum_\nu |S_\nu f|^2)^{1/2}\|_4,\]
which Fefferman proved in \cite{FeffermanIsrael}. At first this appears to run counter to the observation that we can have no equally-spaced forward Littlewood--Paley inequalities for any $p<2$, but on closer inspection it does not violate that principle. Fefferman's argument, which is very simple, exploits the fact that $4$ is an even integer to multiply out the integral implicit on the left hand side to turn it into an $L^2$ expression for which we can use Plancherel's theorem,  and then, crucially, uses {\em the curvature of the circle} to observe that the algebraic sums 
\[ \{{\rm supp} \; \phi_\nu + {\rm supp} \; \phi_{\nu'}\}_{|\nu - \nu'| \lesssim \delta^{-1/2}}\]
have bounded overlap. This, together with arguments for the two-dimensional Kakeya maximal operator and the forward Littlewood--Paley inequalities for equally-spaced intervals when $p \geq 2$ which we discussed above, turns out to be enough to establish the sharp inequality 
\[ \|S^\delta\|_{4-4} \lesssim \log(\delta^{-1})^{1/4}.\]
Let us emphasise that all this is specific to two dimensions, see \cite{CordobaDuke}; the analogous higher-dimensional questions remain open.

Indeed, three of the most important open questions in euclidean harmonic analysis on $\mathbb{R}^n$ for $n \geq 3$ are: (i) $L^p$ bounds for the (Bochner--Riesz) operators $S^\delta$, (ii) $L^p$ bounds for the Kakeya maximal operator (although the corresponding assertion that Kakeya sets have full dimension when $n=3$ has recently been established in spectacular work of Wang and Zahl \cite{WangZahl}), and (iii) $L^p -L^q$ bounds for the Fourier restriction and extension operators for the sphere. These three questions are known to be intimately related, and all of them are resolved when $n=2$, but all are open in all higher dimensions. Moreover, it is known \cite{carberyRLP} that all three would follow from the higher-dimensional analogue
\begin{equation}\label{eq:RLP} \|\sum_\nu S_\nu f\|_{2n/(n-1)} \lesssim  \|(\sum_\nu |S_\nu f|^2)^{1/2}\|_{2n/(n-1)}
\end{equation}
of Fefferman's reverse Littlewood--Paley inequality. This inequality \eqref{eq:RLP}, if true, would therefore be profound. Note that when $n\geq 3$, the index $2n/(n-1)$ is no longer an even integer.

The question of reverse Littlewood--Paley inequalities in the spirit of \eqref{eq:RLP} is equally interesting in related settings. For example, in three dimensions, a $\delta$-neighbourhood of a conical segment $|\xi_3|^2 = |\xi_1|^2 + |\xi_2|^2$, $1 \leq  |\xi_1|^2 + |\xi_2|^2 \leq 2$, can be naturally decomposed into rectangular slabs of dimensions $\delta \times \delta^{1/2} \times 1$, with corresponding smooth multiplier operators $S_\nu$. Once again one has the question of whether   
\[ \|\sum_\nu S_\nu f\|_4 \lesssim  \|(\sum_\nu |S_\nu f|^2)^{1/2}\|_4\]
holds, perhaps up to subpolynomial factors in $\delta^{-1}$. Despite the fact that $4$ is an even integer, this problem had resisted all attempts (see \cite{MSS}) for forty years, and has only recently been resolved in a {\em tour de force} by Guth, Wang and Zhang in \cite{GuthWangZhang}.
The techniques apply not only to the $L^4$-boundedness of the operator $\sum_\nu S_\nu$, but also to local smoothing estimates for dispersive PDE.

Given the difficulty of establishing inequalities such as \eqref{eq:RLP}, it makes sense to consider the corresponding inequalities where the $\ell^2$ and $L^p$ norms are taken in the ``wrong" order, just as with Besov spaces. Thus one is led to study inequalities such as 
\[ \|\sum_\nu S_\nu f\|_p \lesssim  (\sum_\nu \|S_\nu f\|_p^2)^{1/2}\]
(up to sub-polynomial factors in $\delta^{-1}$) for $p \geq 2$. The right-hand side is larger than $  \|(\sum_\nu |S_\nu f|^2)^{1/2}\|_p$ , so there is more scope for establishing positive results, and the range of $p$ for which such might hold is larger than in the classical Littlewood--Paley  setting. Indeed, when $n=2$, the critical $p$ turns out now to be the $L^2$-extension index $6$ instead of $4$. Wolff \cite{Wolff} established sharp estimates of this sort for large $p$ (with the correct power of $\delta$), and that they could be applied to get sharp $L^p$ estimates for radial Fourier multiplier operators for large $p$.
See also \cite{MR2078706} and \cite{MR2664568}. Since then, such estimates, now known as decoupling estimates, have been systematically studied and developed, initially by Bourgain and Demeter \cite{BourgainDemeter} (see also \cite{Demeter}), and have had a huge range of applications, not least of which is the resolution of the Main Conjecture in Vinogradov's Mean Value Theorem, see \cite{BDG}.

\section{Reprise: Almost-everywhere Fourier convergence in $\mathbb{R}^n$ for $n \geq 2$}\label{sec:epilogue}

In this final section we return to the pointwise convergence problems which originally motivated Littlewood and Paley in \cite{LP1}, but now in their higher-dimensional form. 

There are various settings in which we can consider almost-everywhere Fourier convergence -- and the corresponding maximal function -- in dimensions $n \geq 2$. One of the most attractive of these is in euclidean space $\mathbb{R}^n$ where we consider the spherical partial inverse Fourier integrals \eqref{eq:FIRn} for $f$ belonging to various $L^p$ spaces, or perhaps to subspaces consisting of functions with minimal smoothness. For $f \in L^2$ there is no problem in posing the questions since the Fourier transform is an isometric isomorphism of $L^2$.

For $p>2$ there are issues with even making sense of the questions. Indeed, unless $S_R$ maps the Schwartz space $\mathcal{S}(\mathbb{R}^n)$ into $L^{p'}(\mathbb{R}^n)$, it is difficult to see how to define $S_Rf$ even as a tempered distribution for all $f$ belonging to $L^p$. The integral kernel of $S_R$ can be calculated explicitly, and it is in $L^{p'}$ if and only if $p'>2n/(n+1)$. So the question for $L^p$ for $p>2$ only makes sense, even in principle, in the range $p< 2n/(n-1)$, and even then we have to take care. These considerations also underly the occurrence of the index $2n/(n-1)$ in \eqref{eq:RLP}.

For $p<2$, almost-everywhere convergence of $S_{2^k}f$ to $f$ for all $f \in L^p$ would imply via Maurey--Niki\v{s}in--Stein factorisation theory (see  \cite{GCRdeF}) that the maximal operator $f \mapsto \sup_{k \in \mathbb{Z}} |S_{2^k} f|$ would be of weak-type $(p,p)$, and thus $S_1$ would be bounded on some $L^q$ for $q < 2$, which would contradict Fefferman's disc multiplier theorem \cite{Feffermanball}.

What we currently know is summarised by:
\begin{theorem}
Let $2 \leq p < 2n/(n-1)$ and $f \in L^p(\mathbb{R}^n)$. 

(i) 
  Then $S_{2^k}f(x) \rightarrow f(x)$ almost everywhere as $k \to \infty$.

(ii) If for all $\phi \in C^\infty_c$ the function $g = \phi f$ satisfies 
\[ \int_{\mathbb{R}^n} |\widehat{g}(\xi)|^2 (1 + \log^{+} |\xi|)^2 {\rm d} \xi < \infty,\]
then $S_{R}f(x) \rightarrow f(x)$ almost everywhere as $R \to \infty$
\end{theorem}
Part (i) is in \cite{CRdeFV}. The corresponding maximal function estimate for $p=2$ is folklore and can be proved by the same arguments as Theorem~\ref{thm:conv} part (i). Part (ii) is in \cite{CarberySoria} and for $p=2$ is essentially a statement of Rademacher--Men\v{s}ov type. It would be interesting to lower the index $2$ on the logarithmic term to $1$. As a corollary of part (ii) we have almost-everywhere convergence for $f$ in any Sobolev space $L^p_s$ of strictly positive order $s$ when $2 \leq p < 2n/(n-1)$. It is not currently known whether one may recover the result of part (i) when one replaces the lacunary sequence $\{2^k\}_{k \in \mathbb{Z}}$ by a second-order lacunary sequence $\{2^k(1+2^{-j})\}_{k \in \mathbb{Z}, j \in \mathbb{N}}$ (this corresponds to the maximal Hilbert transform in the one-dimensional setting). However, almost-everywhere convergence through a second-order lacunary sequence does hold if $\int_{\mathbb{R}^n} |\widehat{f}(\xi)|^2 (1 + \log^+\log^{+} |\xi|)^2 {\rm d} \xi < \infty$, see \cite{CarbGMT}. 

\bibliographystyle{plain}

\bibliography{LP}
\end{document}